\documentclass[12pt,a4paper,oneside]{amsart}
\usepackage[utf8]{inputenc}
\usepackage{roboto}
\usepackage{anysize}
\usepackage[dvipsnames]{xcolor}
\marginsize{2cm}{2cm}{1.5cm}{1.5cm}
\usepackage{amsfonts, amsmath, amssymb, amsthm, amscd}
\usepackage{upgreek}
\usepackage[english]{babel}
\usepackage{enumitem}
\usepackage{hyperref}
\hypersetup{
	colorlinks   = true, 
	urlcolor     = blue,
	linkcolor    = blue, 
	citecolor    = red 
}
\usepackage[normalem]{ulem}

\numberwithin{equation}{section}
\newtheorem{problem}{Problem}
\newtheorem{thm}{Theorem}[section]
\newtheorem{prp}[thm]{Proposition}
\newtheorem{lem}[thm]{Lemma}
\newtheorem{cor}[thm]{Corollary}

\newtheorem{conj}[thm]{Conjecture}
\theoremstyle{definition}
\newtheorem{dfn}[thm]{Definition}
\newtheorem{rem}[thm]{Remark}

\newcommand{\Href}[2]{\hyperref[#2]{#1~\ref{#2}}}
\newcommand{\Hrefs}[3]{\hyperref[#2]{#1~\ref{#2}} and \hyperref[#3]{\ref{#3}}}
\newcommand{\Hreftitle}[2]{\texorpdfstring{
    \hspace{-4pt}\Href{#1}{#2}}{#1~\ref{#2}}}

\def\phi{\varphi}
\def\epsilon{\varepsilon}
\def\alpha{\upalpha}
\newcommand{\R}{\mathbb{R}}

\newcommand{\braces}[1]{\left\{#1\right\}}
\newcommand{\norm}[1]{\left\|#1\right\|}

\newcommand{\conv}[1]{\operatorname{conv} #1 }

\newcommand{\iprod}[2]{\left\langle#1,#2\right\rangle}

\newcommand{\st}{:\;}

\newcommand{\enorm}[1]{\left|#1\right|}

\DeclareMathOperator{\power}{power}

\newcommand{\noshow}[1]{}

\newcommand{\transpose}[1]{{#1}^{\top}}

\providecommand{\parenth}[1]{\left(#1\right)}
\providecommand{\braces}[1]{\left\{#1\right\}}

\newcommand{\diam}{\operatorname{diam}}
\newcommand{\origin}{\mathbf{0}}

\newcommand{\bigO}[1]{\mathcal{O}\!\parenth{#1}}

\newcommand{\crosssqdist}[2]{d^2\!\parenth{#1, #2}}
\newcommand{\crosssq}[1]{d^2\!\parenth{#1}}
\newcommand{\crosssqdev}[3]{d^2_{#3}\!\parenth{#1, #2}}

\newcommand{\tuplesize}{k}
\newcommand{\colorn}{n}
\newcommand{\circumradname}{\operatorname{Rad}}
\newcommand{\circumrad}[1]{\circumradname\!\parenth{#1}}
\newcommand{\circumradsq}[1]{\circumradname^2\!\parenth{#1}}
\newcommand{\tverbergcrossdiam}{\operatorname{Tv}_{cid}}
\newcommand{\tverbergcol}{\operatorname{Tv}_{c}}
\newcommand{\distransname}{complete system of disjoint transversals}
\newcommand{\completesystem}[1]{\{#1\}}

\DeclareMathOperator{\arcosh}{arcosh}

\title{Tight colorful no-dimensional Tverberg theorem}
\author[P. Barabanshchikova, G. Ivanov, A.~Polyanskii]{{Polina Barabanshchikova and Grigory Ivanov and Alexander~Polyanskii}}

\address{Polina Barabanshchikova,
\newline\hphantom{iii} Department of Computer Science, Aalto University, 02150 Espoo, Finland
}
\email{\href{mailto:polina.barabanshchikova@aalto.fi}{polina.barabanshchikova@aalto.fi}}

\address{Grigory~Ivanov: Pontif\'icia Universidade Cat\'olica do Rio de Janeiro\\
Departamento de Matem\'atica\\
Rua Marqu\^es de S\~ao Vicente, 225\\
Edif\'{\i}cio Cardeal Leme, sala 862\\
22451-900 G\'avea, Rio de Janeiro, Brazil}
\email{\href{mailto:grimivanov@gmail.com}{grimivanov@gmail.com}}

\address{Alexander Polyanskii,
\newline\hphantom{iii} Department of Mathematics, Emory University, Atlanta, GA, 30322, US
}
\email{\href{mailto:apolian@emory.edu}{apolian@emory.edu}}
\urladdr{\url{http://polyanskii.com}}


\begin{document}
\begin{abstract}
We study colorful no-dimensional Tverberg-type problems and obtain several optimal results.
A colorful no-dimensional Tverberg-type theorem provides a bound on a radius $R$ such that, for any pairwise disjoint $k$-element subsets
$Q_1,\dots,Q_n$ of a normed space, there exists a partition of $Q_1\cup\cdots\cup Q_n$ into disjoint transversals $\{P_1,\dots,P_k\}$ for which a ball of radius $R$ intersects the convex hull of each $P_i$ ($1\le i\le k$).
Our methods are deterministic and dimension-free, and they are unified by optimizing two functionals: a quadratic \emph{selection} functional whose local maximizers produce a complete system of disjoint transversals, and a convex \emph{intersection} functional that certifies a common point.

First, in the Euclidean setting we bound $R$ in terms of the Chebyshev radii (minimal enclosing-ball radii) of the color classes $Q_1,\dots,Q_n$.
A key observation is a ``combinatorial'' subadditivity of the squared Chebyshev radius:
given sequences $X=(x_1,\dots,x_k)$ and $Y=(y_1,\dots,y_k)$ of points in a Euclidean space, contained in balls of radii $R_X$ and $R_Y$ (not necessarily with the same center), one can reenumerate $Y$ so that the pointwise-sum sequence $Z=(x_1+y_1,\dots,x_k+y_k)$ is contained in a ball of radius $R_Z$ satisfying
\[
R_Z^2 \le R_X^2 + R_Y^2 .
\]
As a corollary, we obtain the best-possible bound
\[
R \le \frac{1}{\sqrt{2n}}\sqrt{\frac{k-1}{k}}\,
\max_{1\le i\le n} \operatorname{diam}(Q_i).
\]
Our algorithm returns the desired disjoint transversals
 in overall time $\mathcal{O}(nk^3)$.
 
Second, we develop a complementary approach based on the inter-color diameter and extend the framework to obtain no-dimensional colorful Tverberg-type results in the hyperbolic setting and in Banach spaces.
\end{abstract}

\maketitle

\section{Introduction}
\label{section introduction}
In 1966, Helge Tverberg~\cite{tverberg1966} proved his seminal theorem asserting that \textit{for any $(k-1)(d+1)+1$ points in $\mathbb{R}^d$, there is a partition into $k$ parts such that their convex hulls intersect}. This result became a cornerstone for a vast direction in discrete geometry. We refer the reader to the recent survey~\cite{Baranysoberon2018tverberg} by B\'ar\'any and Sober\'on on Tverberg-type results and open problems.

{One of the central open problems in the area is the colorful Tverberg conjecture, which remains open in full generality; see \cite[Section~3]{Baranysoberon2018tverberg}. To state it, we first {recall the notion of} transversals.

For multisets $Q_1,\dots,Q_{\colorn}$, a \emph{transversal} is an $\colorn$-tuple 
$(p_1,\dots,p_{\colorn})$ with $p_i\in Q_i$ for each~$i$. 
{Slightly abusing notation,} we identify $(p_1,\dots,p_{\colorn})$ with its underlying set 
$\{p_1,\dots,p_{\colorn}\}$ when taking convex hulls.
Two transversals are \emph{disjoint} if, for every $i$, they {contain} distinct elements of $Q_i$, that is, they do not share any entry from the same $Q_i$.
If each $Q_i$ has $\tuplesize$ elements, a collection of $\tuplesize$ pairwise disjoint transversals is called a \emph{\distransname} of $Q_1,\dots,Q_{\colorn}$, that is, for {each} $1\le i\le \colorn$, the $\tuplesize$ chosen entries across the transversals coincide with $Q_i$ itself.
}

\begin{conj}[colorful Tverberg conjecture]
\label{conj:colorful_Tverberg}
    For multisets $Q_1, \dots, Q_{d+1}$ of $\tuplesize$ points each in $\mathbb{R}^d$, there is a \distransname \ $\completesystem{P_1,\dots, P_\tuplesize}$ such that the convex hulls of $P_i$ intersect.
\end{conj}

In 2020, Adiprasito \textit{et al.}~\cite{adiprasito2020theorems} proposed a dimension-free variant of the Tverberg problem, in both the uncolored and colorful settings. Their formulation removes the dependence on the ambient dimension from the size of the point set by relaxing the intersection condition: instead of requiring the convex hulls to share a common point, one asks that they all meet a common ball of bounded radius. Their colorful result~\cite[Theorem~6.1]{adiprasito2020theorems} can be formulated as follows.

\begin{prp}[No-dimensional colorful Tverberg theorem]
For multisets $Q_1,\dots,Q_{\colorn}$ of $\tuplesize$ points each in a Euclidean space, there is  a \distransname \  $\completesystem{P_1,\dots, P_\tuplesize}$ and a ball of radius
\begin{equation}\label{eq:ABMT_Euclidean_Tverberg_bound}
\frac{1+\sqrt{2}}{\sqrt{n}} \max\limits_{1 \le i \le \colorn} \diam Q_i
\end{equation}
that intersects the convex hull of each of $P_i,$ $1\le i \le \tuplesize$.
\end{prp}

Their proof uses an averaging argument from the probabilistic method. In particular, the ball can be chosen to be centered at the centroid of the multiset $Q_1\cup\dots\cup Q_{\colorn}$, and it contains the centroids of the transversals. 
{Here and throughout the paper, the centroid of a multiset of vectors is the arithmetic mean of its elements.}
The ``centroid approximation'' viewpoint arises naturally in approximation and sparsification problems addressed via probabilistic techniques; see, for example, Maurey’s lemma~\cite{pisier1980remarques}, Rudelson’s sparsification of sums of positive semidefinite matrices~\cite{rudelson1999random}, and the resolution of the Kadison--Singer problem~\cite{marcus2015interlacing}. As shown by the second author~\cite{Ivanov2021}, combining the approach of~\cite{adiprasito2020theorems} with Maurey’s lemma yields a similar no-dimensional Tverberg-type result in Banach spaces of type $p$, $p>1$. Finally, the bound in~\eqref{eq:ABMT_Euclidean_Tverberg_bound} is optimal up to the constant factor $1+\sqrt{2}$, which cannot be reduced below $1/\sqrt{2}$ uniformly in $\tuplesize$.

To clarify what we mean by the optimal constant, we pose the following no-dimensional colorful Tverberg-type problem in the more general setting of Banach spaces.

\begin{problem}
Given a Banach space $X$ and integers $\tuplesize\ge 2$ and $\colorn\ge 2$, find the smallest value $\tverbergcol=\tverbergcol(X,\tuplesize,\colorn)$
\footnote{The third author used the notation $R(X,\tuplesize,\colorn)$ in the survey~\cite{polyanskii2024} instead of $\tverbergcol(X,\tuplesize,\colorn)$. Since $R$ is standard for a radius, we adopt the notation $\tverbergcol(X,\tuplesize,\colorn)$.}
such that for any $\tuplesize$-point multisets $Q_1,\dots,Q_{\colorn}\subset X$, there exists a \distransname \ $\{P_1,\dots,P_{\tuplesize}\}$ 
of $Q_1,\dots,Q_\colorn$ 
and a ball of radius
\[
\tverbergcol\cdot \max\limits_{1 \le i \le \colorn} \diam Q_i
\]
that intersects the convex hull of each $P_i$, $1\le i \le \tuplesize$.
\end{problem}

\smallskip

Our first main result gives a sharp no-dimensional Tverberg-type bound in Euclidean spaces with the optimal constant.
We suppress {the dependence on} the ambient dimension throughout, since all bounds proved in the paper are no-dimensional. For definiteness, we write $\ell_2$ for the separable infinite-dimensional real Hilbert space.

\begin{thm}\label{thm:nodim_colorful_tverberg_diam_tight}
For $\tuplesize$-point multisets $Q_1,\dots,Q_{\colorn}$ in a Euclidean space, there exists a \distransname \ $\completesystem{P_1,\dots,P_{\tuplesize}}$ of $Q_1,\dots,Q_\colorn$ and a ball of radius
\begin{equation*}
\frac{1}{\sqrt{2}} \sqrt{\frac{k-1}{k}}
\parenth{\frac{1}{\diam^{2} Q_1} + \cdots + \frac{1}{\diam^{2} Q_{\colorn}}}^{-\frac{1}{2}}
\end{equation*}
that intersects the convex hull of each $P_i$, $1\le i \le \tuplesize$.

Moreover, this bound is best possible. There exist multisets $Q_1,\dots,Q_\colorn$ in $\ell_2$ for which no smaller radius can be guaranteed, i.e., for every \distransname , no ball of smaller radius intersects all the convex hulls.

In particular, for any integers $\tuplesize \ge 2$ and $\colorn \ge 2$, we have
\[
\tverbergcol(\ell_2,\tuplesize,\colorn)
= \frac{1}{\sqrt{2\colorn}}\,
\sqrt{\frac{\tuplesize-1}{\tuplesize}}.
\]
\end{thm}

While \Href{Theorem}{thm:nodim_colorful_tverberg_diam_tight} resolves the Euclidean case precisely, in our development it will in fact appear as an almost immediate corollary of more general no-dimensional Tverberg-type theorems proved for a different objective functional. Namely, instead of the diameter, we use the minimal radius of a ball containing a set (also called the \emph{Chebyshev radius} or \emph{circumradius}). We formulate two optimal colorful no-dimensional Tverberg-type statements in \Href{Subsection}{subsec:optimality_cicrumradius_results}, one of which concerns ``centroid approximation'' and provides an optimal bound on the radius of the ball enclosing the centroids of the transversals. 
Next, we {state} the key ingredient of the proof, the ``combinatorial'' subadditivity of the squared Chebyshev radius.
\begin{thm}\label{thm:combinatorial_subadd_circumradius}
Let $X=(x_1,\ldots, x_\tuplesize)$ and $Y=(y_1,\ldots,y_\tuplesize)$ be sequences of points in a Euclidean space contained in balls of respective radii $R_X$ and $R_Y$. It is possible to reenumerate the points of $Y$ so that the sequence $Z=(x_1+y_1,\ldots, x_\tuplesize + y_\tuplesize)$ is contained in a ball of radius $R_Z$ with
\[
R_Z^2 \leq {R_X^2  + R_{Y}^2}.
\]
\end{thm}

We note three points.
First, instead of the maximum diameter, our bounds use the suitable average, thereby strengthening the results of \cite{adiprasito2020theorems}.
Second, the proof is non-probabilistic: the desired \distransname \  is a ``local'' maximizer of a quadratic functional.
The idea of proving Tverberg-type results using quadratic cost functionals goes back to two proofs of Tverberg’s theorem -- one by Tverberg and Vrećica~\cite{Tverberg1993} and another by Roudneff~\cite{roudneff2001partitions}. Our proof, however, builds on a substantially reworked version of this idea.
In addition to optimizing a quadratic functional, we also optimize an auxiliary function to locate a ``Tverberg point.''\footnote{In the classical setting this means a common point of convex hulls, while in our setting it refers to the center of a ball intersecting the convex hulls, or a common point of balls.} This refinement was introduced in \cite[Theorem~4]{Pirahmad2022}. We also point the reader to papers~\cite{barabanshchikova2024intersecting, barabanshchikova2024intersectingballs} of the {first} and third authors, where a similar two-objective technique was applied to other Tverberg-type problems. For a broader overview of no-dimensional Tverberg-type problems, we refer to the recent survey~\cite{polyanskii2024}.
{Third, the desired \distransname \ can be found in time polynomial in $\colorn\tuplesize$, thus improving upon the results of \cite{Choudhary2022}}.\smallskip

The drawback of \Href{Theorem}{thm:nodim_colorful_tverberg_diam_tight}  is that it is purely Euclidean. Conceptually, we seek the center of a smallest enclosing ball (the \emph{Chebyshev center}), but an issue arises: in normed spaces of dimension at least three, the property that the Chebyshev center of every compact set lies in its convex hull characterizes inner-product spaces; it holds if and only if the space is Hilbert (see the Garkavi--Klee characterization \cite{Garkavi1964,klee1960circumspheres}).\smallskip

For this reason, we would like to introduce a new variation of the colorful no-dimensional Tverberg
problem by changing the objective functional. Specifically, we replace the maximum diameter with an \emph{inter-color diameter}, defined for $\colorn \ge 2$ finite subsets $Q_1,\dots,Q_{\colorn}$ of a Banach space $X$ by
\[
\diam(Q_1,\dots,Q_{\colorn})
:=
\max\{\norm{x-y}_X \;:\; x\in Q_i,\ y\in Q_j \text{ with } 1 \leq i < j \leq \colorn \}.
\]
Notably, this notion of diameter also appears in a colorful variant of Jung’s theorem studied by Akopyan~\cite{Akopyan2013}. He showed that for any $\colorn$ sets $Q_1,\dots,Q_{\colorn}$ in $\ell_2$, there exists a ball of radius $\diam(Q_1,\dots,Q_{\colorn})/\sqrt{2}$ that covers all sets except one. Clearly, this ball also intersects every transversal.

We now state the following no-dimensional colorful Tverberg-type problem:
\begin{problem}
    \label{problem colorful}
    Given a Banach space $X$ and integers $ \tuplesize \ge 2$ and $\colorn\ge 2$, find the smallest value $\tverbergcrossdiam=\tverbergcrossdiam(X, \tuplesize,\colorn)$ such that for any $\tuplesize$-point multisets $Q_1,\dots,Q_{\colorn}\subset X$, there exists a \distransname \ $\completesystem{P_1,\dots,P_{\tuplesize}}$ of $Q_1,\dots,Q_\colorn$ and a ball of radius
    \[
        \tverbergcrossdiam\cdot \diam(Q_1,\dots,Q_{\colorn})
    \]
    that intersects the convex hull of each $P_i$.
\end{problem}

By the triangle inequality, for any configuration of points, we have
\[
\max_{1\le j\le \colorn} \diam Q_j \;\le\; 2\,\diam(Q_1,\dots,Q_{\colorn}).
\]
Clearly, $\diam(Q_1,\dots,Q_{\colorn})$ cannot be bounded above by
$\max_{1\le j\le \colorn} \diam Q_j$.

We obtain the best possible bound for $\tverbergcrossdiam$ in the Euclidean case.
\begin{thm}
\label{thm:cross-diam_Tverberg}
For any integers $\tuplesize \ge 2$ and $\colorn \ge 2$,
\(
\tverbergcrossdiam(\ell_2,\tuplesize,\colorn)=\frac{1}{\sqrt{2\colorn}}.
\)
\end{thm}
{\Href{Theorem}{thm:cross-diam_Tverberg} bounds the radius in terms of the inter-color diameter $\diam(Q_1,\dots,Q_{\colorn})$, whereas \Href{Theorem}{thm:nodim_colorful_tverberg_diam_tight} uses $\max_{1\le i\le \colorn}\diam Q_i$. Since $\diam(Q_1,\dots,Q_{\colorn})$ can exceed $\max_{1\le i\le \colorn}\diam Q_i$ by an arbitrarily large factor, \Href{Theorem}{thm:cross-diam_Tverberg} does not strengthen \Href{Theorem}{thm:nodim_colorful_tverberg_diam_tight}. Nevertheless, it remains of independent interest.}

First, using an additional shifting trick, we obtain a version of \Href{Theorem}{thm:nodim_colorful_tverberg_diam_tight} with a slightly weaker multiplicative constant; see \Href{Subsection}{subsec:inter_diam_results}.
Second, we believe that the approach used to prove \Href{Theorem}{thm:cross-diam_Tverberg} applies to other norms, most notably the $L_p$ norms and the Schatten norms. For example, in \Href{Section}{sec:hyperbolic_space} we derive a no-dimensional colorful Tverberg-type result in hyperbolic space by suitably modifying the proof. {Third, our proof essentially follows the approach in the proof of \Href{Theorem}{thm:nodim_colorful_tverberg_diam_tight} and optimizes two functionals: one indicates a desired complete system of transversals, while another --- a ``Tverberg point''.}

We conclude the introduction with an open problem, which can be viewed as a relaxed version of the colorful Tverberg conjecture.
\begin{problem}[weak colorful Tverberg problem]
Determine a non-trivial upper bound on $R$ such that the following holds. For $\tuplesize$-point multisets $Q_1,\dots,Q_{d+1}\subset \R^d$, there  exists {a \distransname \ $\completesystem{P_1,\dots,P_{\tuplesize}}$ of $Q_1,\dots,Q_{d + 1}$} and a ball of radius
\[
R\cdot \diam P
\]
that intersects the convex hull of each $P_i,$ $1\le i \le \tuplesize$.
\end{problem}

The rest of the paper is organized as follows.
In the next section, we outline the geometric ideas underlying our results in Euclidean space.
In \Href{Section}{sec:circumradius_tverberg_results}, we prove \Href{Corollary}{cor:optimal_circum_rad_bound} and \Href{Corollary}{cor:tverberg_for_centroids} and establish the sharpness of these bounds together with \Href{Theorem}{thm:nodim_colorful_tverberg_diam_tight}. These theorems follow from a more delicate result on the convexity properties of the minimum radius of a ball intersecting all transversals.
In \Href{Section}{sec:inter-diameter results}, we present an alternative approach via inter-diameter bounds and derive \Href{Theorem}{thm:cross-diam_Tverberg} together with several corollaries. We then extend this approach to the hyperbolic setting and prove a hyperbolic analogue of \Href{Theorem}{thm:cross-diam_Tverberg} in \Href{Section}{sec:hyperbolic_space}. Finally, in \Href{Section}{sec:tverberg_banach_spaces} we discuss the Banach-space case and pose several open problems.

\subsection*{Notation}
{For a finite multiset $Q=\{q_1,\dots,q_m\}$ in a vector space, the centroid of $Q$, denoted by $c_Q$, is defined as the average of its elements, that is,}
\[
c_Q = \frac{q_1+\cdots+q_m}{m}.
\]

When the ambient space is Euclidean, the Euclidean norm of a vector $a$ is denoted by $\enorm{a}$.
We work only with finite multisets, and any auxiliary points we introduce will lie in the linear hull (even in the convex hull) of the corresponding multisets. Thus, in the Euclidean parts one may, without loss of generality, regard the ambient space as $\ell_2$.

{By default, the sum $\sum_{1\le i\ne j\le t} a_{ij}$ runs over $t(t-1)$ \emph{ordered} pairs, while $\sum_{1\le i<j\le t} a_{ij}$ runs over $\binom{t}{2}$ \emph{unordered} pairs.}

\section{General idea behind the Euclidean results}
\label{sec:idea_Euclidean_no-dim_Tverberg}
We now delve into {the ideas behind} the proofs of our Euclidean results, namely \Href{Theorem}{thm:cross-diam_Tverberg} and \Href{Theorem}{thm:combinatorial_subadd_circumradius}.
As mentioned in the Introduction, we use two functionals: the first selects the desired \distransname , while the second shows that a suitable ball intersects the convex hull of each transversal.

\medskip
\noindent\textbf{The selection functional.}
All the disjoint transversals $\completesystem{P_1,\dots,P_{\tuplesize}}$ we seek will be (local) maximizers of quadratic functionals in the initial data $(Q_1,\dots,Q_{\colorn})$.
Writing $P_i=\{x_{i1},\dots,x_{i\colorn}\}$ with $x_{ij}\in Q_j$ the unique point chosen by $P_i$ from color $j$, a typical quadratic functional has the form
\[
\operatorname{F}_1(P_1,\dots,P_{\tuplesize})
=\sum_{1\le i,j\le \tuplesize}\ \sum_{1\le p,q\le \colorn} \mu_{ijpq}\,\iprod{x_{ip}}{x_{jq}} \;+\; C,
\]
where the coefficients $\mu_{ijpq}$ are fixed and $C=C(Q_1,\dots,Q_{\colorn})$ is a suitable constant.

In practice, we introduce these functionals via a geometric description that implicitly absorbs the additive constant, using
\[
2\,\iprod{a}{b}=\enorm{a}^2+\enorm{b}^2-\enorm{a-b}^2
\]
to rewrite $\operatorname{F}_1$ as a sum of squared lengths. 
For example, in the two-color case ($\colorn=2$) relevant to \Href{Theorem}{thm:combinatorial_subadd_circumradius}, we take
\[
{\operatorname{F}_1(P_1,\dots,P_{\tuplesize})
=\enorm{ x_{11} - x_{12}}^2+\cdots+\enorm{ x_{\tuplesize1}-x_{\tuplesize2}}^2.}
\]

Throughout, it suffices that $P_1,\dots,P_{\tuplesize}$ are \emph{local} maximizers of the corresponding functional, meaning that swapping $x_{it}$ and $x_{jt}$ for any fixed $t$ does not increase its value.

\medskip
\noindent\textbf{Why a maximization helps a minimization.}
One might wonder how maximizing a functional helps with problems that are, at heart, 
minimizations -- after all, we seek balls of minimal radius. 
Informally, the exchange-optimality of the selection functional yields a family of pairwise \emph{swap inequalities}, and these inequalities are exactly what we need when we pass to the second, \emph{Tverberg-point functional}.

\medskip
\noindent\textbf{The Tverberg-point functional.}
The choice of the Tverberg-point functional is more delicate. In fact, once the appropriate Tverberg-point functional is identified, the corresponding ``selection'' functional follows almost immediately. For \Href{Theorem}{thm:cross-diam_Tverberg}, we study the intersection of the Euclidean balls
\[
B_i=\braces{\,x\in \ell_2 \st \enorm{x-c_i}^2\le r_i^2\,},\qquad 1\le i\le \tuplesize.
\]
We consider
\[
\operatorname{F}_2(x):=\max_{1\le i\le \tuplesize}\,\parenth{\enorm{x-c_i}^2-r_i^2}.
\]
The balls $B_1,\dots,B_{\tuplesize}$ have a common point if and only if $\inf_x \operatorname{F}_2(x)\le 0$; equivalently, it suffices to show that the minimum value of $\operatorname{F}_2$ is nonpositive. Since $\operatorname{F}_2$ is convex (a pointwise maximum of convex quadratics), standard subdifferential calculus yields a necessary optimality condition at a minimizer $z$,
\[
\origin \in \partial \operatorname{F}_2(z),
\]
which can be rewritten as a balancing condition. This condition, in turn, transforms into the family of swap inequalities that are exactly those satisfied by a local maximizer of the selection functional.

{It is worth noting that in \Href{Theorem}{thm:combinatorial_subadd_circumradius} and its corollaries, the Tverberg-point functional is used implicitly: for a given set of points, it assigns to each point the maximum distance to the set, and its minimizer is precisely the center of the minimal enclosing ball.}\bigskip

\section{Chebyshev radius results}
\label{sec:circumradius_tverberg_results}
We start by proving \Href{Theorem}{thm:combinatorial_subadd_circumradius}.
We recall related definitions and results.

\begin{dfn}
For a subset $A$ of a Euclidean space $\mathcal{H}$, the \emph{Chebyshev radius} $\circumrad{A}$ is the minimal radius of a ball containing $A$, that is,
\[
\circumrad{A}:=\inf_{c\in \mathcal{H}}\ \sup_{x\in A}\,\enorm{x-c}.
\]
A point $c$ such that the ball of radius $\circumrad{A}$ centered at $c$ contains $A$ is called the \emph{Chebyshev center} of $A$.
\end{dfn}

The existence and uniqueness of the Chebyshev center (and, hence, of the minimal enclosing ball) is classical \cite[Theorem~2]{Garkavi1964}:

\begin{lem}\label{lem:chebyshev_center}
Every compact set $A$ in a Euclidean space has a unique Chebyshev center, and it lies in the convex hull of $A$.
\end{lem}

We shall also use the following optimization characterization due to Alexander \cite{alexander1984circumdisk}.

\begin{prp}[{\cite[Theorem~2]{alexander1984circumdisk}}]
\label{prp:alexander_circumrad}
Let $A=\{a_1,\dots,a_m\}$ be a finite subset of a Euclidean space. Its Chebyshev radius satisfies
\[
2\circumradsq{A}=\max\ \transpose{\lambda}\,M\,\lambda,
\]
where the maximum is over all $\lambda\in\R^m$ with nonnegative coordinates summing to one, and $M$ is the matrix with entries $M[i,j]=\enorm{a_i-a_j}^2$.
\end{prp}

\begin{proof}[Proof of \Href{Theorem}{thm:combinatorial_subadd_circumradius}]
Relabel the indices of $Y$ so that
\[
\sum\limits_{1 \le i \le \tuplesize} \enorm{x_i - y_i}^{2}
\]
is maximized. Swapping $y_i$ and $y_j$ for any $i\neq j$, the maximality yields
\[
\enorm{x_i - y_i}^{2} + \enorm{x_j - y_j}^{2}
\ \ge\ 
\enorm{x_i - y_j}^{2} + \enorm{x_j - y_i}^{2}.
\]
Expanding the squares, this is equivalent to
\[
2\iprod{x_i}{y_i} + 2 \iprod{x_j}{y_j}\ \le\ 
2 \iprod{x_i}{y_j} + 2\iprod{x_j}{y_i}
\ \Longleftrightarrow\
2 \iprod{x_i-x_j}{\,y_i-y_j}\ \le\ 0.
\]
Therefore, adding $\enorm{x_i-x_j}^{2} + \enorm{y_i-y_j}^{2}$, we obtain
\begin{equation}
\label{eq:trans}
\enorm{(x_i+y_i)-(x_j+y_j)}^{2}
\ \le\ \enorm{x_i-x_j}^{2} + \enorm{y_i-y_j}^{2}.
\end{equation}

Let $Z=\{z_1,\dots,z_\tuplesize\}$ with $z_i:=x_i+y_i$. Denote by $M_Z$ the $\tuplesize\times\tuplesize$ matrix with entries $M_Z [i,j]=\enorm{z_i-z_j}^2$, similarly define $M_X$ and $M_Y$.
The inequality \eqref{eq:trans} translates to inequality 
\[
M_Z [i,j] \leq M_X [i,j] + M_Y [i,j].
\]
Hence, by \Href{Proposition}{prp:alexander_circumrad},
\[
\circumradsq{Z} \le
\circumradsq{X} + \circumradsq{Y} 
\ \le\ R_X^2+R_Y^2,
\]
{which finishes the proof.}
\end{proof}

\begin{rem}
{Roman} Karasev pointed out that, after \eqref{eq:trans}, the proof of \Href{Theorem}{thm:combinatorial_subadd_circumradius} can be completed via Kirszbraun’s extension theorem \cite{kirszbraun1934zusammenziehende}.

Let \(E\) be the Euclidean space under consideration, and set
\(P=(x_1\oplus y_1,\dots,x_{\tuplesize}\oplus y_{\tuplesize})\subset E\oplus E\).
By the Pythagorean theorem, this set is contained in a ball of radius
\(R_Z=\sqrt{R_X^{2}+R_Y^{2}}\).

Define
\[
f:\ P\to Z,\qquad f(x_i\oplus y_i)=x_i+y_i.
\]
Inequality \eqref{eq:trans} shows that \(f\) is \(1\)-Lipschitz. By Kirszbraun’s theorem \cite{kirszbraun1934zusammenziehende}, \(f\) extends to a \(1\)-Lipschitz map on all of \(E\oplus E\); in particular it is defined at the center \(c\) of the enclosing ball of \(P\). Therefore, for each \(i\),
\[
\enorm{f(x_i\oplus y_i)-f(c)}\ \le\ \enorm{(x_i\oplus y_i)-c}\ \le\ R_Z,
\]
so \(Z=f(P)\) is contained in the ball of radius \(R_Z\) centered at \(f(c)\).
\end{rem}

By the standard induction argument, we obtain
\begin{thm}
    \label{thm:comb_subadd_general}
  Let $X_i=(x_{1i},\ldots, x_{\tuplesize i}),$ $1 \le i \le \colorn,$ be sequences of points in a Euclidean space. It is possible to reenumerate the points of $X_j,$ $2 \le j \le \colorn$ so that the sequence 
  \[
  Z=(x_{11}  + \dots + x_{1 \colorn},\ldots, x_{\tuplesize 1}   
  + \dots + x_{ \tuplesize \colorn} ) 
  \]
  satisfies the inequality
\[
\circumradsq{Z} \le \circumradsq{X_1} + \dots + \circumradsq{X_\colorn}.
\]
\end{thm}

\subsection{Corollaries of combinatorial subadditivity}
\label{subsec:optimality_cicrumradius_results}
In this subsection, we derive several optimal no-dimensional Tverberg-type results, including \Href{Theorem}{thm:nodim_colorful_tverberg_diam_tight}. They follow from \Href{Theorem}{thm:comb_subadd_general}.

For $\tuplesize$-point multisets $Q_1,\dots,Q_{\colorn}$, set
\begin{equation}
\label{eq:rad_centroid_def}
      R_c = \frac{\big(\circumradsq{Q_1}+\dots+\circumradsq{Q_{\colorn}}\big)^{1/2}}{\colorn}
\end{equation}
and
\begin{equation}
\label{eq:rad_hm_def}
R_{hm} =
    \parenth{\frac{1}{\circumradsq{Q_1}} + \cdots +
    \frac{1}{\circumradsq{Q_{\colorn}}}}^{-\frac{1}{2}}.
\end{equation}
Geometrically, these quantities give the optimal radii in our no-dimensional Tverberg-type bounds. We begin with examples illustrating their appearance.

As a direct consequence of \Href{Lemma}{lem:chebyshev_center}, we obtain the following folklore fact.
\begin{lem}\label{lem:chebyshev_center_charachterization}
{Let $Q$ be a finite multiset in $\ell_2$} contained in a ball $B$ of radius $R$ centered at $c$. The following assertions are equivalent:
\begin{enumerate}
  \item $R$ is the Chebyshev radius $\circumradname{Q}$ of $Q$, and $c$ is the Chebyshev center of $Q$.
  \item There exist points $q_1,\dots,q_t$ of $Q$ such that $\enorm{q_i-c}=R$ for all $1\le i\le t$, and
  $c=\sum_{1\le i\le t}\theta_i q_i$ for some nonnegative coefficients $\theta_1,\dots,\theta_t$ with $\sum_{1\le i\le t}\theta_i=1$.
\end{enumerate}
\end{lem}

\begin{lem}\label{lem:optimality_nodim_Tverberg}
Let $Q_1,\dots,Q_{\colorn}\subset \ell_2$ be the vertex sets of regular simplices of dimension $\tuplesize-1$, whose affine hulls are pairwise orthogonal and whose centroids are at the origin. Then the following holds for any \distransname\ $\completesystem{P_1,\dots,P_{\tuplesize}}$ of $Q_1,\dots,Q_{\colorn}$:
\begin{enumerate}
  \item\label{item:lem_optimality_tverberg_Ch_centroid}
  The Chebyshev radius of the set of centroids $\{c_{P_1},\dots,c_{P_{\tuplesize}}\}$ equals $R_c$ given by \eqref{eq:rad_centroid_def}.
  \item\label{item:lem_optimality_tverberg_Ch_rad}
  There is no ball of radius strictly less than $R_{hm}$ given by \eqref{eq:rad_hm_def} that intersects all the convex hulls of $P_1,\dots,P_{\tuplesize}$.
\end{enumerate}
\end{lem}

\begin{proof}
The first statement is immediate: for any transversal $P$ of $Q_1,\dots,Q_{\colorn}$, its centroid $c_P$ satisfies
\(
\enorm{c_P}=R_c.
\)
Moreover, the centroid of the multiset $\{c_{P_1},\dots,c_{P_{\tuplesize}}\}$, for any $\tuplesize$ pairwise disjoint transversals $\completesystem{P_1,\dots,P_{\tuplesize}}$ of $Q_1,\dots,Q_{\colorn}$, is the origin. The claim follows from \Href{Lemma}{lem:chebyshev_center_charachterization}.

For the second assertion, let $P=\{q_1,\dots,q_{\colorn}\}$ be any transversal with $q_i\in Q_i$ for $1\le i\le \colorn$. The distance from the origin to the affine hull of $P$ is exactly $R_{hm}$. Moreover, the orthogonal projection $h$ of the origin onto the affine hull of $P$ lies in $\conv P$ and can be written as $\alpha_1 q_1+\dots+\alpha_{\colorn} q_{\colorn}$. By invariance under orthogonal transformations, the coefficient $\alpha_i$ depends only on the index $i$ and not on the particular choice of $q_i\in Q_i$. Hence, the sum of the orthogonal projections of the origin onto the affine hulls of any $\tuplesize$ pairwise disjoint transversals $\completesystem{P_1,\dots,P_{\tuplesize}}$ is zero.

Consequently, for any point $x\in\ell_2$ there exists an index $i$ with $1\le i\le \tuplesize$ such that the distance from $x$ to the affine hull of $P_i$ is at least $R_{hm}$. The claim follows.
\end{proof}

As direct corollaries of \Href{Theorem}{thm:comb_subadd_general} and the identity
$\circumradsq{\lambda A} = \lambda^2 \circumradsq{A}$, we obtain that $R_c$ and $R_{hm}$ are the optimal radii in the corresponding Tverberg-type problems for the Chebyshev radius.

\begin{cor}\label{cor:optimal_circum_rad_bound}
For $\tuplesize$-point multisets $Q_1,\dots,Q_{\colorn}$ in a Euclidean space, there exists a \distransname\ $\completesystem{P_1,\dots,P_{\tuplesize}}$ of $Q_1,\dots,Q_\colorn$ and a ball of radius $R_{hm}$ given by \eqref{eq:rad_hm_def}
that intersects the convex hull of each $P_i$, $1\le i\le \tuplesize$.

Moreover, this bound is best possible: there exist multisets $Q_1,\dots,Q_\colorn$ in $\ell_2$ for which no smaller radius can be guaranteed, i.e., for every \distransname, no ball of smaller radius intersects all the convex hulls.
\end{cor}

\begin{proof}
Set
\[
\lambda_i =
\frac{R_{hm}^2}{\circumradsq{Q_i}},
\quad 1 \le i \le \colorn,
\]
and apply \Href{Theorem}{thm:comb_subadd_general} to the sets 
$\lambda_i Q_i,$ $1 \le i \le \colorn.$ This yields the stated bound since 
$\sum_{1 \le i \le \colorn}\lambda_i = 1$.
Optimality follows from \eqref{item:lem_optimality_tverberg_Ch_rad} of \Href{Lemma}{lem:optimality_nodim_Tverberg}.
\end{proof}

Since the ``centroid approximation'' is of independent interest, we state the optimal result below.

\begin{cor}
\label{cor:tverberg_for_centroids}
For $\tuplesize$-point multisets $Q_1,\dots,Q_{\colorn}$ in a Euclidean space, there exists a \distransname\ $\completesystem{P_1,\dots,P_{\tuplesize}}$ of $Q_1,\dots,Q_\colorn$ and a ball of radius $R_{c}$ given by \eqref{eq:rad_centroid_def}.

Moreover, this bound is best possible: there exist multisets $Q_1,\dots,Q_\colorn$ in $\ell_2$ for which no smaller radius can be guaranteed, i.e., for every \distransname, no ball of smaller radius contains all the centroids.
\end{cor}

\begin{proof}
Choose $\lambda_1=\cdots=\lambda_{\colorn}=\frac{1}{\colorn}$ and apply \Href{Theorem}{thm:comb_subadd_general} to the sets 
$\lambda_i Q_i,$ $1 \le i \le \colorn.$ This yields the stated bound.
Optimality follows from \eqref{item:lem_optimality_tverberg_Ch_centroid} of \Href{Lemma}{lem:optimality_nodim_Tverberg}.
\end{proof}

By classical Jung’s theorem \cite{jung1901}, a $\tuplesize$-point multiset of diameter $D$ in a Euclidean space can be covered by a ball of radius
\begin{equation}
\label{eq:Jung_inequality}
R \le D\,\sqrt{\frac{k-1}{2k}}.
\end{equation}
Applying this within \Href{Corollary}{cor:optimal_circum_rad_bound}, we recover the optimal constant in \eqref{eq:ABMT_Euclidean_Tverberg_bound}.

\begin{proof}[Proof of \Href{Theorem}{thm:nodim_colorful_tverberg_diam_tight}]
The upper bound follows from \Href{Corollary}{cor:optimal_circum_rad_bound} together with Jung’s inequality \eqref{eq:Jung_inequality}. Optimality, as well as the identity
\[
\tverbergcol(\ell_2,\tuplesize,\colorn)=\frac{1}{\sqrt{2\colorn}}\sqrt{\frac{\tuplesize-1}{\tuplesize}},
\]
follow directly from the example in \Href{Lemma}{lem:optimality_nodim_Tverberg} and the relation
\[
\circumrad{\Delta}=\sqrt{\frac{\tuplesize-1}{2\tuplesize}}\,\diam \Delta
\]
for a regular simplex $\Delta$ with $\tuplesize$ vertices.
\end{proof}

\subsection{Algorithmic aspects of the proof}

Since we work with finitely many points, we may assume, without loss of generality, that
$Q_1,\dots,Q_{\colorn}$ lie in the finite-dimensional subspace given by the linear hull of
$Q_1\cup\dots\cup Q_{\colorn}$. Under this reduction, we can find a desired
\distransname\ $\completesystem{P_1,\dots,P_{\tuplesize}}$ and a witnessing ball in polynomial time.

\medskip
In the proof of \Href{Theorem}{thm:nodim_colorful_tverberg_diam_tight}, we apply \Href{Theorem}{thm:comb_subadd_general}; more precisely, we apply the induction step, that is,  \Href{Theorem}{thm:combinatorial_subadd_circumradius}, $\colorn-1$ times to the corresponding multisets.
\Href{Theorem}{thm:combinatorial_subadd_circumradius} returns a perfect matching (equivalently, permutations $X$ and $Y$) maximizing
\[
\sum_{1 \le i \le \tuplesize}\enorm{\lambda x_i-(1-\lambda)y_i}^{2}
\]
for the relevant $\lambda \in [0,1]$.
This is a linear assignment problem with cost matrix
$C_{ij}=\enorm{\lambda x_i-(1-\lambda)y_j}^{2}$, solvable by the Hungarian algorithm in
$\bigO{\tuplesize^3}$ time. Since the base step is invoked $\colorn-1$ times in the
inductive construction, we obtain $\tuplesize$ disjoint transversals
$\completesystem{P_1,\dots,P_{\tuplesize}}$ in overall time $\bigO{\colorn\tuplesize^3}$.

\medskip
Computing Chebyshev centers is a classical optimization problem. Using Alexander’s distance-matrix formulation mentioned above, one can compute the Chebyshev center for a given $\tuplesize$-point set in time $\bigO{\tuplesize^3}$.

\section{Inter-diameter results}
\label{sec:inter-diameter results}
The goal of this section is to prove \Href{Theorem}{thm:no-dim_Tverberg_cross_diam} and 
\Href{Corollary}{cor:tverberg_diam}. We start by introducing yet another distance-like functional.

For a set $P=\{x_1,\dots,x_{\colorn}\}$, let $\crosssq{P}$ and $B_2(P)$ be a constant and a ball defined as follows
\[
    \crosssq{P}:=\frac{1}{2\colorn^2(\colorn-1)}\sum_{1 \le i \neq j \le \colorn}\enorm{x_i-x_j}^2
\]
and
\[
    B_2(P):=\braces{\,x\in \ell_2:\ \enorm{x-c_P}^2\le \crosssq{P}\,}.
\]

Consider two sets $A=\{a_1,\dots,a_{\colorn}\}$ and $B=\{b_1,\dots,b_{\colorn}\}$. 
For each $1\le m\le \colorn$, let us introduce the notation
\[
\crosssqdev{A}{B}{m}
=
\sum_{\substack{1\le t\le \colorn}, t \ne m}
\parenth{\enorm{a_t-b_m}^2+\enorm{a_m-b_t}^2}
\]
and
\[
\crosssqdist{A}{B}
=
\frac{1}{4\colorn^2(\colorn-1)}
\sum_{1\le m\le \colorn} d_m^2(A,B).
\]
One can easily see that $\crosssqdist{A}{A}=\crosssq{A}$.

Our main theorem -- which implies the corresponding statements from the introduction -- is as follows.
\begin{thm}\label{thm:no-dim_Tverberg_cross_diam}
For $\tuplesize$-point multisets $Q_1,\dots,Q_{\colorn}\subset \ell_2$, there exist $\tuplesize$ pairwise disjoint transversals $\completesystem{P_1,\dots,P_{\tuplesize}}$ such that the balls $B_2(P_i)$ share a common point.
\end{thm}

\begin{proof}
 The argument proceeds in two steps. First, we show that there are disjoint transversals $\completesystem{P_1,\dots,P_{\tuplesize}}$ satisfying 
\begin{equation}
    \label{eq: main inequality}
    \crosssq{P_i}+\crosssq{P_j}-2\crosssqdist{P_i}{P_j} \ge 0
    \quad \text{for any } 1\le i\ne j\le \tuplesize.
\end{equation}

To prove the first claim, we choose disjoint transversals $\completesystem{P_1,\dots,P_{\tuplesize}}$ maximizing the sum
\begin{equation}
    \label{equation:maxsum}
    \crosssq{P_1} + \dots + \crosssq{P_{\tuplesize}}. 
\end{equation}
Suppose that \eqref{eq: main inequality} does not hold for some $1\le i\ne j\le \tuplesize$. Then there is an index $1\le t\le \colorn$ such that
\[      
    \crosssqdev{P_i}{P_i}{t} + \crosssqdev{P_j}{P_j}{t} - 2\crosssqdev{P_i}{P_j}{t} < 0.
\]

Denote by $x_{it}$ and $x_{jt}$ the elements of the singletons $P_i\cap Q_t$ and $P_j\cap Q_t$, respectively.  Interchanging these between $P_i$ and $P_j$, we obtain new transversals 
\[
   P_i' = 
    \bigl(P_i\setminus \{x_{it}\}\bigr) \cup \{x_{jt}\}
    \quad\text{and}\quad 
   P_j' = 
    \bigl(P_j\setminus \{x_{jt}\}\bigr) \cup \{x_{it}\}.
\]
Therefore,
\[
    \crosssqdev{P_i'}{P_i'}{t} + \crosssqdev{P_j'}{P_j'}{t} 
    = 2 \crosssqdev{P_i}{P_j}{t} 
    > \crosssqdev{P_i}{P_i}{t} + \crosssqdev{P_j}{P_j}{t},
\]
and thus,
\[
    \crosssq{P_i'} + \crosssq{P_j'} > \crosssq{P_i} + \crosssq{P_j},
\]
a contradiction with the maximality of~\eqref{equation:maxsum}. Henceforth, we assume that disjoint transversals $\completesystem{P_1,\dots,P_{\tuplesize}}$ satisfy~\eqref{eq: main inequality}.

\medskip

For a point $x$ and a transversal $P$, denote by $\power(x,P)$ the \textit{power} of $x$ with respect to the ball $B_2(P)$, that is,
\[
    \power(x,P):= 
    \enorm{x-c_{P}}^2 - \crosssqdist{P}{P}.
\]
Consider the function $f:\ell_2\to \R$ defined by
\[
    f(x)=\max_{1\le i\le \tuplesize} \big\{\power(x,P_i)\big\}.
\]
If this function attains a non-positive value at some point $x$, then the point 
$x$ is a common point of the balls $B_2(P_i)$, $ 1 \le i \le \tuplesize$. Suppose that these balls do not share a common point. Then, since $f$ is convex, it attains a positive value at its minimum point. Without loss of generality, assume that this minimum point is the origin $\origin$. 

For simplicity, assume that $f(\origin)=\power(\origin,P_i)$ if and only if $1\le i\le l$, where $l\le \tuplesize$, and $P_i=\{x_{i1},\dots,x_{i\colorn}\}$, where $x_{it}$ is the only element of the singleton $P_i\cap Q_t$. A routine computation shows that, for each $1\le i\le l$,
\begin{equation}
\label{equation positive sum}
0<\frac{2\colorn^2(\colorn-1)}{2\colorn-1}\, f(\origin)
= \frac{2\colorn^2(\colorn-1)}{2\colorn-1}\,\power(\origin,P_i)
= \sum_{\substack{1\le t\ne m\le \colorn}} 
\iprod{x_{it}}{x_{im}} .
\end{equation}

Since $f$ is a pointwise maximum of convex functions, its subdifferential {$\partial f$ at $\mathbf 0$}, defined as
\[
\partial f(z)=\{y\in \ell_2:\ \iprod{y}{x-z} \le f(x)-f(z) \text{ for all }x\in \ell^2\},
\]
coincides with the convex hull of the gradients 
\[
    \nabla \power (\origin,P_i)=-2c_{P_i}, \quad \text{for all } 1\le i\le l,
\]
namely those indices for which the maximum is achieved at $\origin$.
Since $\origin$ is the minimum point of $f$, we have $\origin \in \partial f(\origin)$, that is,
\[
    \origin=\colorn\sum_{1\le i\le l} \lambda_i c_{P_i}
     =\sum_{1\le i\le l} \lambda_i\sum_{1\le t\le \colorn} x_{it}
     =\sum_{1\le t\le \colorn}\ \sum_{1\le i\le l}\lambda_i x_{it},
\]
where $\lambda_i\ge 0$ and $\sum\limits_{1\le i\le l} \lambda_i=1$. 
Thus,
\[
0=
\enorm{\sum_{1\le t\le \colorn}\ \sum_{1\le i\le l}\lambda_i x_{it}}^2
=
\sum_{1\le t\le \colorn}\enorm{\sum_{1\le i\le l}\lambda_i x_{it}}^2
+\sum_{\substack{1\le t\ne m \le \colorn}}
\iprod{\sum_{1\le i\le l}\lambda_i x_{it}}{\sum_{1\le j\le l}\lambda_j x_{jm}}.
\]

Hence,
\begin{align*}
    0\ge {}&\sum_{1\le t\ne m\le \colorn}
             \iprod{\sum_{1\le i\le l}\lambda_i x_{it}}{\sum_{1\le j\le l}\lambda_j x_{jm}}\\
    ={}& \sum_{1\le t\ne m\le \colorn} \sum_{1\le i \le l} \lambda_i^2 \iprod{x_{it}}{x_{im}}
        + \sum_{1\le t\ne m\le \colorn} \sum_{1\le i\ne j\le l} \lambda_i \lambda_j \iprod{x_{it}}{x_{jm}}\\
    ={}&\sum_{1\le i \le l} \lambda_i^2\sum_{1\le t\ne m\le \colorn}\iprod{x_{it}}{x_{im}}
        +\sum_{1\le i\ne j\le l} \lambda_i \lambda_j \sum_{1\le t\ne m \le \colorn} \iprod{x_{it}}{x_{jm}}.
\end{align*}

By~\eqref{equation positive sum}, we get
\[
    \sum_{1\le i\ne j\le l} \lambda_i \lambda_j \sum_{1\le t\ne m \le \colorn} \iprod{x_{it}}{x_{jm}}<0.
\]
This implies that there exist indices $1\le i\ne j\le l$ such that 
\[
    \sum_{1\le t\ne m\le \colorn} \big(\iprod{x_{it}}{x_{jm}}+\iprod{x_{it}}{x_{jm}}\big)
    <0
    <
    \sum_{1\le t\ne m\le \colorn} \big(\iprod{x_{it}}{x_{im}} + \iprod{x_{jn}}{x_{jm}}\big),
\]
where the last inequality follows from \eqref{equation positive sum}. Therefore,
\[
   2 \crosssqdist{P_i}{P_j}
   > 
   \crosssqdist{P_i}{P_i} +
   \crosssqdist{P_j}{P_j},
\]
a contradiction.
\end{proof}

\begin{rem}
As was mentioned in \Href{Section}{sec:idea_Euclidean_no-dim_Tverberg}, any \distransname \ $\completesystem{P_1,\dots, P_k}$ \emph{locally} maximizing the sum 
    \begin{equation}
    \label{equation invariant}
     \crosssq{P_1} +\dots + \crosssq{P_\tuplesize}
    \end{equation}
    satisfy this statement. Here ``locally maximizing'' means that there is no swap between two of the transversals $P_i$ and $P_j$ that increases the sum.
\end{rem}

\subsection{Proofs of  \Hreftitle{Theorem}{thm:cross-diam_Tverberg} and related results}
\label{subsec:inter_diam_results}
\begin{proof}[Proof of \Href{Theorem}{thm:cross-diam_Tverberg}]
The result follows from \Href{Theorem}{thm:no-dim_Tverberg_cross_diam} and the following upper bound on the radius of $B_2(P_i)$:
\[
   \crosssq{P_i} \le \frac{\colorn(\colorn-1)\,\diam^2 P_i}{2\colorn^2(\colorn-1)}
    \le \frac{\diam^2(Q_1,\dots,Q_{\colorn})}{2\colorn},
\]
which completes the proof.
\end{proof}

Now, let us show the trick allowing us to derive  \Href{Theorem}{thm:nodim_colorful_tverberg_diam_tight} alike result from 
\Href{Theorem}{thm:no-dim_Tverberg_cross_diam}. The idea is that we can shift all the sets.
For two $\colorn$-tuples $A = (a_1, \dots, a_\colorn)$ and $B =  (b_1, \dots, b_\colorn)$ in a linear space,  we define their difference $A - B$ as $(a_1 - b_1, \dots, a_\colorn - b_\colorn).$ 
\begin{lem}\label{lem:shifting_trick_tverberg_cross_diam}
For $\tuplesize$-point multisets $Q_1,\dots,Q_{\colorn}\subset \ell_2$ and  an $\colorn$-tuple 
$T = (t_1, \dots, t_{\colorn})$ of vectors of $\ell_2$, there exist $\tuplesize$ disjoint transversals $\completesystem{P_1,\dots,P_{\tuplesize}}$ such that the translated balls $B_2(P_i - T)+c_T$ (centered at $c_{P_i}$) share a common point.
\end{lem}
\begin{proof}
Using \Href{Theorem}{thm:no-dim_Tverberg_cross_diam} for the shifted sets $Q_1-t_1,\dots,Q_{\colorn}-t_{\colorn}$, we obtain that the balls 
\[
B_2(P_i - T)=\bigl\{x\in \ell_2:\ \enorm{x-c_{P_i - T}}^2 \le  \crosssq{P_i - T}\bigr\}
\]
share a common point, denoted by $o$. Since $c_{P-T}=c_P-c_T$, the point $o+c_T$ is a common point of the balls
\[
\bigl\{x\in \ell_2:\ \enorm{x-c_{P_i}}^2 \le \crosssq{P_i - T}\bigr\},
\quad 1 \le i \le \tuplesize.
\]
\end{proof}
\begin{cor}\label{cor:tverberg_diam}
For $\tuplesize$-point multisets $Q_1,\dots,Q_{\colorn}$ in a Euclidean space, there exist $\tuplesize$ disjoint transversals $\completesystem{P_1,\dots,P_{\tuplesize}}$ and a ball of radius
\begin{equation*}
\sqrt{2}\,\frac{\sqrt{\diam^2 Q_1+\cdots+\diam^2 Q_{\colorn}}}{\colorn}
\end{equation*}
{that intersects the convex hull of each of $P_i,$ $1\le i \le \tuplesize$.}
\end{cor}
\begin{proof}

The result follows from \Href{Lemma}{lem:shifting_trick_tverberg_cross_diam} and the following upper bound on the radius of $B_2(P - T)$ {for any transversal $P = (x_1, \dots, x_\colorn)$} in the case $T = (c_{Q_1}, \dots, c_{Q_\colorn})$:
\[
   \crosssq{P {- T}} \le \frac{1}{2\colorn^2(\colorn-1)} 
   \sum_{1 \le i \neq j \le \colorn}\enorm{(x_i - c_{Q_i}) - (x_j - c_{Q_j})}^2 
    \le
\]
\[
\frac{1}{2\colorn^2(\colorn-1)} 
   \sum_{1 \le i \neq j \le \colorn} 2(\diam^2 Q_i + \diam^2 Q_j) = 2 
   \frac{\diam^2 Q_1+\cdots+\diam^2 Q_{\colorn}}{\colorn^2}.
\]

\end{proof}

\section{A hyperbolic variation of the no-dimensional Tverberg theorem}
\label{sec:hyperbolic_space}

In $\R\times \ell_2$, consider the symmetric bilinear form $p$ defined by
\[
p(x,y):=-x_0y_0+\iprod{x_\infty}{y_\infty},
\]
for vectors $x=(x_0,x_\infty)$ and $y=(y_0,y_\infty)$. The associated quadratic form is $q(x):=p(x,x)$.

Define the hyperboloid
\[
\mathbb H=\{x=(x_0,x_\infty)\in \R\times \ell_2:\ q(x)=-1 \text{ and } x_0>0\}.
\]
Equipped with the hyperbolic distance
\[
d(x,y)=\arcosh\!\big(-p(x,y)\big),\qquad x,y\in\mathbb H,
\]
the space $(\mathbb H,d)$ is the hyperboloid model of the infinite-dimensional hyperbolic space.

For a $t$-point set $X\subset \mathbb H$ with $t\ge 2$, define the hyperbolic ball induced by $X$ as
\[
B_H(X):=\bigl\{\,z\in \mathbb H:\ \sum_{x\ne y\in X} p(x-z,y-z)\le 0\,\bigr\}.
\]
Let $c_X$ be the (ambient) centroid of $X$ in $\R\times \ell_2$. Then $B_H(X)$ is a hyperbolic ball centered at $c_X/\sqrt{\enorm{q(c_X)}}\in \mathbb H$. Indeed, for $z\in\mathbb H$,
\[
\sum_{x\ne y\in X} p(x-z,y-z)
=\sum_{x\ne y\in X} p(x,y)-2(t-1)\sum_{x\in X} p(z,x)-t(t-1),
\]
so there exists $r\in \R$ such that $z\in B_H(X)$ if and only if $p(z,c_X)\ge r$, equivalently
\[
d\!\left(z,\ \frac{c_X}{\sqrt{\enorm{q(c_X)}}}\right)\le \arcosh\!\left(-\frac{r}{\sqrt{\enorm{q(c_X)}}}\right).
\]
In particular, when $t=2$, this yields the hyperbolic ball centered at the midpoint of the geodesic between $x_1,x_2\in X$ with radius $d(x_1,x_2)/2$.

We now present the no-dimensional colorful Tverberg theorem in hyperbolic space.
\begin{thm}
\label{theorem hyperbolic}
For $\tuplesize$-point sets $Q_1,\dots,Q_{\colorn}$ in $\mathbb H$, there exist \distransname \ $\completesystem{P_1,\dots,P_{\tuplesize}}$ such that the balls $B_H(P_i)$ share a common point.
\end{thm}
\begin{proof}
    Consider two transversals $A=\{a_1,\dots,a_{\colorn}\}$ and $B=\{b_1,\dots,b_{\colorn}\}$, where $a_i,b_i\in Q_i$. Define
\[
p(A,B):=\sum_{\substack{1\le i\ne j\le \colorn}} p(a_i,b_j).
\]

Choose disjoint transversals $\completesystem{P_1,\dots,P_{\tuplesize}}$ minimizing
\[
p(P_1,P_1)+\dots+p(P_{\tuplesize},P_{\tuplesize}).
\]
Following the same argument as in the proof of \Href{Theorem}{thm:no-dim_Tverberg_cross_diam} for~\eqref{eq: main inequality}, with $d^2$ replaced by $-p$, we obtain
\begin{equation}
\label{eq : hyperboloid transversals criterion}
p(P_i,P_i)+p(P_j,P_j)-2p(P_i,P_j)\le 0
\quad \text{for any } 1\le i\ne j\le \tuplesize.
\end{equation}

It remains to show that $P_1,\dots,P_{\tuplesize}$ satisfy the conditions of the theorem.

For any $z\in \ell_2$, set $z_H = \big(\sqrt{1+\iprod{z}{z}},\,z\big)\in \mathbb H$ and define
\[
f:\ell_2\to \R,\qquad
f(z):=\max_{1\le i\le \tuplesize} f_i(z),
\quad \text{where}\quad
f_i(z):=\sum_{x\ne y\in P_i} p(x-z_H,y-z_H).
\]
Using the bilinearity of $p$,
\[
f_i(z)=p(P_i,P_i)-2(\colorn-1)\sum_{x\in P_i} p(x,z_H)-\colorn(\colorn-1)
\ \ge\ p(P_i,P_i)-\colorn(\colorn-1),
\]
where the inequality follows from the Cauchy--Schwarz-type observation below.

\begin{prp}\label{proposition:time-like}
Let $x = (x_0,x_{\infty})$ and $y = (y_0,y_{\infty}) \in \R \times \ell_2$ with $x_0,y_0>0$ and $q(x),q(y)<0$. Then:
\begin{enumerate}[label=(\roman*)]
\item $p(x,y)<0$;
\item $p(x,y)^2\ge q(x)\,q(y)$, with equality if and only if $x = \lambda y$ for some 
$\lambda \in \R$.
\end{enumerate}
\end{prp}

Thus, each $f_i$ is bounded below. Since $z\mapsto \sqrt{1+\iprod{z}{z}}$ is convex, each $f_i$ is convex, and hence so is $f$ (as a pointwise maximum). Therefore, $f$ attains its minimum at some $z\in \ell_2$. If $f(z)\le 0$, then $z_H\in B_H(P_i)$ for all $i$, and we are done. Suppose instead that $f(z)>0$.

Without loss of generality, there exists $l\le \tuplesize$ such that $f(z)=f_i(z)$ if and only if $1\le i\le l$. In particular, for each $1\le i\le l$,
\begin{equation}
\label{eq : hyperboloid balanced set}
0<f(z)=f_i(z)=\sum_{\substack{x,y\in P_i}} p(x-z_H,y-z_H).
\end{equation}
Since $z$ minimizes $f$, we have $\origin\in \partial f(z)$, so there exist $\lambda_i\ge 0$ with $\sum\limits_{1 \le i \le l} \lambda_i = 1$ such that
\begin{equation}
\label{eq: hyperbolic balancing condition}
\origin=\sum\limits_{1 \le i \le l} \lambda_i \nabla f_i(z).
\end{equation}

Write $P_i=\{x_{i1},\dots,x_{i\colorn}\}$, where $x_{ij}$ is the unique element of $P_i\cap Q_j$. We use the following technical inequality (proved in \Href{Section}{sec:proof_technical_ineq_hyperboloid}):
\begin{equation}
\label{eq:technical_ineq_hyperboloid}
\sum_{1\le m\ne t\le \colorn}
p\!\left(\sum\limits_{1 \le i \le l} \lambda_i (x_{it}-z_H),\ \sum\limits_{1 \le i \le l} \lambda_i (x_{im}-z_H)\right)\ \le\ 0.
\end{equation}

From \eqref{eq:technical_ineq_hyperboloid} we get
\begin{align*}
0 &\ge
\sum_{1\le m\ne t\le \colorn}
p\!\left(\sum\limits_{1 \le i \le l} \lambda_i (x_{it}-z_H),\ \sum\limits_{1 \le i \le l} \lambda_i (x_{im}-z_H)\right)\\
&=
\sum\limits_{1\le m\ne t\le \colorn}
\sum\limits_{1 \le i \le l} \lambda_i^2\, p(x_{it}-z_H, x_{im}-z_H)
+\sum_{1\le m\ne t\le \colorn}\sum_{\substack{1\le i\ne j\le l}} \lambda_i\lambda_j\, p(x_{it}-z_H, x_{jm}-z_H)\\
&=
\sum\limits_{1 \le i \le l} \lambda_i^2 \sum_{1\le m\ne t\le \colorn} p(x_{it}-z_H, x_{im}-z_H)
+\sum_{\substack{1\le i\ne j\le l}} \lambda_i\lambda_j \sum_{1\le m\ne t\le \colorn} p(x_{it}-z_H, x_{jm}-z_H).
\end{align*}
By \eqref{eq : hyperboloid balanced set},
\[
\sum_{\substack{1\le i\ne j\le l}} \lambda_i\lambda_j \sum_{1\le m\ne t\le \colorn} p(x_{it}-z_H, x_{jm}-z_H)\ <\ 0.
\]
Hence there exist $1\le i\ne j\le l$ such that
\[
\sum_{1\le m\ne t\le \colorn} 2\,p(x_{it}-z_H, x_{jm}-z_H)
<
\sum_{1\le m\ne t\le \colorn} p(x_{it}-z_H, x_{im}-z_H)
+\sum_{1\le m\ne t\le \colorn} p(x_{jt}-z_H, x_{jm}-z_H),
\]
where the right-hand inequality follows from \eqref{eq : hyperboloid balanced set}. Using bilinearity of $p$ and canceling the terms involving $z_H$ on both sides, we obtain
\[
2\sum_{1\le m\ne t\le \colorn} p(x_{it},x_{jm})
<
\sum_{1\le m\ne t\le \colorn} p(x_{it},x_{im})
+\sum_{1\le m\ne t\le \colorn} p(x_{jt},x_{jm}),
\]
and therefore
\[
2\,p(P_i,P_j) < p(P_i,P_i)+p(P_j,P_j),
\]
contradicting \eqref{eq : hyperboloid transversals criterion}.
\end{proof}
\subsection{\texorpdfstring{Proof of \eqref{eq:technical_ineq_hyperboloid}}
{Proof of the technical inequality}}
\label{sec:proof_technical_ineq_hyperboloid}
First, we derive a closed-form expression for $z_H$. 
Recall that $z_H = \parenth{\sqrt{1+\iprod{z}{z}}, \,z}$. The gradient of $f_i$ at $z$ is
\[
\nabla f_i(z)
= -2(\colorn-1)\sum_{(x_0,x_{\infty})\in P_i}
\left(x_{\infty}-x_0\,\frac{z}{\sqrt{1+\iprod{z}{z}}}\right).
\]
Using the balancing condition \eqref{eq: hyperbolic balancing condition}, we find
\[
z
= c^{-1}\sum_{1\le i\le l}\lambda_i\!\!\sum_{(x_0,x_{\infty})\in P_i}\!\! x_{\infty}
\quad\text{and}\quad
\sqrt{1+\iprod{z}{z}}
= c^{-1}\sum_{1\le i\le l}\lambda_i\!\!\sum_{(x_0,x_{\infty})\in P_i}\!\! x_0,
\]
where $c\in\R_{+}$ is determined by $q(z_H)=-1$. Denoting $s_t=\sum_{1\le i\le l}\lambda_i\,x_{it}$, we obtain
\[
z_H = c^{-1}\sum_{1\le t\le \colorn} s_t,
\qquad
c=\sqrt{-\,q\!\left(\sum_{1\le t\le \colorn} s_t\right)}.
\]

With this notation, the left-hand side of \eqref{eq:technical_ineq_hyperboloid} becomes
\begin{align*}
\sum_{1\le m\ne t\le \colorn} p(s_t - z_H, s_m - z_H)
&= \sum_{1\le m\ne t\le \colorn} p(s_t,s_m)
   -2(\colorn-1)\,p\!\left(z_H,\sum_{1\le t\le \colorn}s_t\right)
   +\colorn(\colorn-1)\,q(z_H)\\
&= q\!\left(\sum_{1\le t\le \colorn}s_t\right)
   -\sum_{1\le t\le \colorn} q(s_t)
   +2(\colorn-1)\sqrt{-\,q\!\left(\sum_{1\le t\le \colorn}s_t\right)}
   -\colorn(\colorn-1)\\
&= -\Bigg(\sqrt{-\,q\!\left(\sum_{1\le t\le \colorn}s_t\right)}-(\colorn-1)\Bigg)^2
   -\sum_{1\le t\le \colorn} q(s_t) - (\colorn-1).
\end{align*}

By Proposition~\ref{proposition:time-like}, for any $x,y\in\mathbb H$ we have $p(x,y)\le -1$. Hence, for each $1\le t\le \colorn$,
\[
q(s_t)=\sum_{1\le i,j\le l}\lambda_i\lambda_j\,p(x_{it},x_{jt})
\;\le\; -\Big(\sum_{1\le i\le l}\lambda_i\Big)^2=-1.
\]
Let $\beta_t:=\sqrt{-\,q(s_t)}\ge 1$. Proposition~\ref{proposition:time-like} also gives
$p(s_t,s_m)^2\ge q(s_t)\,q(s_m)=\beta_t^2\beta_m^2$, and therefore
\[
-\,q\!\left(\sum_{1\le t\le \colorn}s_t\right)
\;\ge\; \sum_{1\le m\ne t\le \colorn}\beta_t\beta_m+\sum_{1\le t\le \colorn}\beta_t^2
= \Big(\sum_{1\le t\le \colorn}\beta_t\Big)^2
\;\ge\; \colorn^2.
\]

Consequently,
\begin{align*}
\sum_{1\le m\ne t\le \colorn} p(s_t - z_H, s_m - z_H)
&\le
-\Big(\sum_{1\le t\le \colorn}\beta_t-(\colorn-1)\Big)^2
+\sum_{1\le t\le \colorn}\beta_t^2 - (\colorn-1)\\
&=
-\Big(\sum_{1\le t\le \colorn}(\beta_t-1)\Big)^2
+\sum_{1\le t\le \colorn}(\beta_t-1)^2
\;\le\; 0,
\end{align*}
which completes the proof of \eqref{eq:technical_ineq_hyperboloid}.
\subsection{Spherical case}
Note that there is no similar result in spherical geometry. 
Since a sphere has positive curvature, two balls induced by the opposite edges of a spherical square do not intersect. For instance, let $Q_1$ and $Q_2$ be pairs of antipodal points lying at the intersection of the unit sphere with the first and the second axes, respectively. For any partition of these points into disjoint transversals $P_1$ and $P_2$, the spherical segments connecting points in $P_1$ and $P_2$ have antipodal midpoints, and each segment has a length of $\pi/2$. Consequently, two spherical circles with diameters given by these segments do not intersect.

\section{No-dimensional colorful Tverberg-type results in Banach spaces}
\label{sec:tverberg_banach_spaces}
We believe that the results of \Href{Section}{sec:inter-diameter results} extend to certain Banach spaces after modifying the distance-like functional $\crosssq{\cdot}$ as follows.

Let $X$ be a Banach space with norm $\norm{\cdot}_X$, and fix $q\in[1,2]$. For a transversal $P$, define the constant $\Delta^q(P)$ and the ball $B_X(P,q)$ by
\[
    \Delta^q(P):=\frac{1}{2\colorn^q(\colorn-1)}\,
    {\displaystyle\sum_{x\ne y\in P}\norm{x-y}_X^{\,q}},
    \qquad
    B_X(P,q):=\braces{\,x\in X:\ \norm{x-c_P}_X^{\,q}\le \Delta^q(P)\,}.
\]
Note that when $q=2$ and $X=\ell_2$, we recover $\Delta^q(P)=\crosssq{P}$. We conjecture that \Href{Theorem}{thm:no-dim_Tverberg_cross_diam} extends to $L_p$ spaces, $p\in[1,\infty)$, via a multicomponent Riesz--Thorin interpolation argument (see \cite[Chapter~4]{Wells1975}). Moreover, in view of the derandomization of Maurey’s lemma in \cite{ivanov2021approximate}, it is natural to study the functional $\Delta^q$ in uniformly smooth Banach spaces.

We now present our main result for point sets in an arbitrary Banach space $X$.
The proof is a straightforward modification of the argument of \Href{Theorem}{thm:no-dim_Tverberg_cross_diam} in the simplest case $q=1$.
To lighten notation, we henceforth write $\Delta$ for $\Delta^{1}$ and $B_X(P)$ for $B_X(P,1)$.

\begin{thm}
\label{thm:banach_tverberg_constant_bound}
Let $X$ be a Banach space equipped with the norm $\norm{\cdot}_X$. For $\tuplesize$-point multisets $Q_1,\dots,Q_{\colorn}\subset X$, there exists a 
\distransname \ $\completesystem{P_1,\dots, P_\tuplesize}$ of $Q_1,\dots,Q_{\colorn}$
such that the balls $B_X(P_i)$ are pairwise intersecting.
\end{thm}
\begin{proof}
    Let us show that disjoint transversals $\completesystem{P_1,\dots,P_{\tuplesize}}$ maximizing
\begin{equation}
\label{equation:maxPi}
\Delta(P_1)+\dots+\Delta(P_{\tuplesize})
\end{equation}
satisfy the conclusion of the theorem.

We prove that the balls $B_X(P_i)$ and $B_X(P_j)$ intersect for any $1\le i\ne j\le \tuplesize$. Without loss of generality, assume $i=1$, $j=2$, and write
$P_1=\{a_1,\dots,a_{\colorn}\}$, $P_2=\{b_1,\dots,b_{\colorn}\}$ with $a_t,b_t\in Q_t$.

For each $1\le t\le \colorn,$ we have
\begin{equation}
\label{equation:maximality}
\sum_{1\le i\ne j\le \colorn}\!\big(\norm{a_i-b_j}_X+\norm{a_j-b_i}_X\big)
\;\le\;
\sum_{1\le i\ne j\le \colorn}\!\big(\norm{a_i-a_j}_X+\norm{b_i-b_j}_X\big).
\end{equation}
Otherwise, swapping $a_t$ and $b_t$ between $P_1$ and $P_2$ yields new transversals
\[
P_1'=\bigl(P_1\setminus\{a_t\}\bigr)\cup\{b_t\},
\qquad
P_2'=\bigl(P_2\setminus\{b_t\}\bigr)\cup\{a_t\},
\]
for which
\[
\Delta(P_1)+\Delta(P_2)\;<\;\Delta(P_1')+\Delta(P_2'),
\]
contradicting the maximality of \eqref{equation:maxPi}.

Using the triangle inequality and \eqref{equation:maximality}, we bound the distance between the centers $c_{P_1}$ and $c_{P_2}$ of $B_X(P_1)$ and $B_X(P_2)$:
\begin{align*}
\norm{c_{P_1}-c_{P_2}}_X
&=\frac{1}{\colorn} \norm{\sum_{1 \le i \le \colorn}(a_i-b_i)}_X\\
&=\frac{1}{2\colorn(\colorn-1)}
 \norm{\sum_{1\le i\ne j\le \colorn}(a_i-b_j)
        +\sum_{1\le i\ne j\le \colorn}(a_j-b_i)}_X\\
&\le \frac{1}{2\colorn(\colorn-1)}
   \Big(\sum_{1\le i\ne j\le \colorn}\norm{a_i-b_j}_X
       +\sum_{1\le i\ne j\le \colorn}\norm{a_j-b_i}_X\Big)\\
&\le \frac{1}{2\colorn(\colorn-1)}
   \sum_{1\le i\ne j\le \colorn}\norm{a_i-a_j}_X
   +\frac{1}{2\colorn(\colorn-1)}
   \sum_{1\le i\ne j\le \colorn}\norm{b_j-b_i}_X\\
&={\Delta(P_1)+\Delta(P_2)}.
\end{align*}
The right-hand side is the sum of the radii of $B_X(P_1)$ and $B_X(P_2)$, hence these two balls intersect.
\end{proof}
Similarly to the shifting trick in \Href{Lemma}{lem:shifting_trick_tverberg_cross_diam}, we may shift the multisets and obtain the following consequence.
\begin{cor}
Let $X$ be a Banach space equipped with the norm $\norm{\cdot}_X$. For $\tuplesize$-point multisets $Q_1,\dots,Q_{\colorn}\subset X$, there exists a 
\distransname \ $\completesystem{P_1,\dots, P_\tuplesize}$ of $Q_1,\dots,Q_{\colorn}$
such that the diameter of the multiset of centroids
$\braces{c_{P_1},\dots,c_{P_\tuplesize}}$ is at most
\[
2\,\frac{\tuplesize-1}{\tuplesize}\cdot \frac{\diam Q_1+\cdots+\diam Q_{\colorn}}{\colorn}.
\]
\end{cor}
\begin{proof}
Using \Href{Theorem}{thm:banach_tverberg_constant_bound} for the shifted sets 
$Q_1-c_{Q_1},\dots,Q_{\colorn}-c_{Q_\colorn}$, we obtain {$\tuplesize$ disjoint transversals $\completesystem{P_1,\dots,P_{\tuplesize}}$ of $Q_1,\dots,Q_{\colorn}$ such that for any $1\le i\ne j\le \tuplesize$},
\[
\norm{c_{P_i}-c_{P_j}}_X \;\le\; 
\Delta(P_i - C) + \Delta(P_j - C),
\]
where $C=\parenth{c_{Q_1},\dots,c_{Q_\colorn}}$.

On the other hand, for a transversal 
$P=\parenth{x_1,\dots,x_{\colorn}}$ with $x_i\in Q_i$, we have
\[
\begin{aligned}
\Delta(P - C)
&\le
\frac{1}{2\colorn(\colorn-1)} 
\sum\limits_{1 \le i \neq j \le \colorn}
\norm{(x_i - c_{Q_i}) - (x_j - c_{Q_j})}_X \\[2pt]
&\le
\frac{1}{2\colorn(\colorn-1)} 
\sum\limits_{1 \le i \neq j \le \colorn}
\parenth{\norm{x_i - c_{Q_i}}_X + \norm{x_j - c_{Q_j}}_X}
=
\frac{1}{\colorn}\sum\limits_{1 \le i \le \colorn}\norm{x_i - c_{Q_i}}_X.
\end{aligned}
\]
Finally, by the centroid identity, {we have}
\[
x_i - c_{Q_i}
= \frac{1}{\tuplesize}\sum\limits_{y\in Q_i} (x_i - y),
\]
so by the triangle inequality, {we obtain}
\[
\norm{x_i - c_{Q_i}}_X
\le \frac{1}{\tuplesize}\sum\limits_{y\in Q_i}\norm{x_i-y}_X
\le \frac{\tuplesize-1}{\tuplesize}\,\diam Q_i.
\]
Combining the estimates yields
\[
\norm{c_{P_i}-c_{P_j}}_X
\le
2\,\frac{\tuplesize-1}{\tuplesize}\cdot \frac{\diam Q_1+\cdots+\diam Q_{\colorn}}{\colorn},
\]
and the stated bound on the diameter of $\braces{c_{P_1},\dots,c_{P_\tuplesize}}$ follows.
\end{proof}

\subsection{Tverberg property of Banach spaces}
In view of the recent paper \cite{artstein2025b}, it is natural to formulate the following conjecture.

Recall that a Banach space $X$ is said to be of \emph{Rademacher type $p$} if there exists a constant $C>0$ such that for every finite set $\{x_1,\dots,x_m\}\subset X$,
\[
\frac{1}{2^m}
\sum_{(\alpha_1,\dots,\alpha_m)\in\braces{-1,1}^m}
\norm{\sum\limits_{1 \le i \le m}  \alpha_i x_i}_X
\;\le\;
C\left(\sum\limits_{1 \le i \le m} \norm{x_i}_X^{\,p}\right)^{\!1/p}.
\]

We say that a Banach space $X$ has the \emph{inter-diameter Tverberg property} if
\[
\tverbergcrossdiam(X,\tuplesize,\colorn)\to 0
\quad\text{as }\colorn\to\infty,
\]
for each fixed $\tuplesize\ge 2$.

Analogously, we say that a Banach space $X$ has the \emph{Tverberg property} if
\[
\tverbergcol(X,\tuplesize,\colorn)\to 0
\quad\text{as }\colorn\to\infty,
\]
for each fixed $\tuplesize\ge 2$.

\begin{conj}
The following assertions are equivalent for a Banach space $X$:
\begin{enumerate}
    \item $X$ has the Tverberg property.
    \item $X$ has the inter-diameter Tverberg property.
    \item $X$ is of Rademacher type $p$ for some $p>1$.
\end{enumerate}
\end{conj}

To support the conjecture, we compute $\tverbergcrossdiam(\ell_\infty,\tuplesize,\colorn)$ for the space $\ell_\infty$, which is not of Rademacher type $p$ with $p>1$.

\begin{thm}
\label{theorem l^infty}
For any integers $\tuplesize\ge 2$ and $\colorn\ge 2$, we have
\[
\tverbergcrossdiam(\ell_\infty,\tuplesize,\colorn)=\frac{1}{2}.
\]
\end{thm}
\begin{proof}
It is well known that finitely many balls in $\ell_\infty$ have a common point if and only if they are pairwise intersecting. Thus the bound
\[
\tverbergcrossdiam(\ell_\infty,\tuplesize,\colorn)\ \le\ \frac12
\]
follows from \Href{Theorem}{thm:banach_tverberg_constant_bound}, since
\[
    {\Delta(P_i)}\ \le\ \frac{\colorn(\colorn-1)\,\diam P_i}{2\colorn(\colorn-1)}
    \ \le\ \frac{\diam(Q_1,\dots,Q_{\colorn})}{2}.
\]

We now show the reverse inequality
\[
\tverbergcrossdiam(\ell_\infty,\tuplesize,\colorn)\ \ge\ \frac12.
\]
Our argument uses the example from~\cite[last paragraph of Section~1]{Ivanov2021}, which gives the lower bound $R_c(\ell_1,\tuplesize,\colorn)\ge \tfrac12$. We embed this example into $\ell_\infty$.

Let $T:\ell_1\to \ell_\infty$ be an isometric embedding of $\ell_1$ into $\ell_\infty$ (distance-preserving on its image), whose existence is shown, e.g., in~\cite[Theorem~2.5.7]{albiac2006topics}. Define  $\tuplesize$-point sets 
\[
Q_i=\braces{\,T(e_{\tuplesize i+1}),\dots,T(e_{\tuplesize(i+1)})\,},\qquad i=0,1,\dots,\colorn-1,
\]
where $e_j$ denotes the $j$-th standard basis vector of $\ell_1$.
Suppose there exist $\tuplesize$ disjoint transversals $\completesystem{P_1,\dots,P_{\tuplesize}}$ and a ball of radius $R$ intersecting the convex hull of each $P_i$. A direct computation shows that the distance between any point of $\conv P_1$ and any point of $\conv P_2$ is exactly $2$. By the triangle inequality, the distance between these convex hulls is at most $2R$, hence $R\ge 1$. Since $\diam(Q_1,\dots,Q_{\colorn})=2$, we obtain
\[
\tverbergcrossdiam(\ell_\infty,\tuplesize,\colorn)\ \ge\ \frac{R}{\diam(Q_1,\dots,Q_{\colorn})}\ \ge\ \frac12,
\]
which proves sharpness.
\end{proof}

\section*{Acknowledgment}

The authors  thank Alexey Vasileuski for carefully reading the first version of the manuscript, {and Roman Karasev for his in-depth discussion of the paper's results.} 

A.P. is partially supported by the NSF grant DMS 2349045. 
G. I. is supported by Projeto Paz and by CAPES (Coordena\c{c}\~ao de Aperfei\c{c}oamento de Pessoal de N\'ivel Superior) - Brasil, grant number 23038.015548/2016-06.
\bibliographystyle{alpha}
\bibliography{biblio}
\end{document}